\documentclass[11pt,a4paper,leqno,oneside]{article}

\usepackage{amsmath}
\usepackage{amsfonts}
\usepackage{amssymb}
\usepackage{amsthm}
\usepackage{makeidx}
\usepackage{graphicx}
\usepackage{pb-diagram}
\usepackage{epstopdf}
\usepackage{fancyhdr}
\usepackage{enumitem}
\usepackage[all]{xy}

\newtheorem{teo}{Theorem}
\newtheorem{prop}[teo]{Proposition} 
\newtheorem{lem}[teo]{Lemma}
\newtheorem{cor}[teo]{Corollary}
\theoremstyle{definition}
\newtheorem{defin}[teo]{Definition}

\newtheorem{prob}[teo]{Problem}

   {\begin{verse}%
     \small }%
   {\end{verse}}

\title{CYCLIC HILBERT SPACES AND CONNES' EMBEDDING PROBLEM}

\author{Valerio CAPRARO and Florin R\u ADULESCU}

\begin{document}

\maketitle

$\bf{Abstract}$. Let $M$ be a $II_1$-factor with trace $\tau$, the
linear subspaces of $L^2(M,\tau)$ are not just common Hilbert
spaces, but they have additional structure. We introduce the notion
of a cyclic linear space by taking those properties as axioms. In
Sec.2 we formulate the following problem: ''does every cyclic
Hilbert space embed into $L^2(M,\tau)$, for some $M$?''. An
affirmative answer would imply the existence of an algorithm to
check Connes' embedding Conjecture. In Sec.3 we make a first step
towards the answer of the previous question.

\tableofcontents

\section{Cyclic Hilbert spaces}
Let $M$ be a finite factor with unique normalized trace $\tau$ and
let $L^2(M,\tau)$ be the Hilbert space obtained by taking the
closure of the vector space $M$ with respect to the inner product
$(x,y)=\tau(y^*x)^{\frac{1}{2}}$. Consider a finite-dimensional real
Hilbert subspace $H\subseteq M_{sa}\subseteq L^2(M,\tau)$ containing
the identity. Observe that $H$ is not just a common Hilbert space,
but it has additional structure.
\begin{prop}\label{properties}
The mapping $\ll a\otimes b,c\otimes d\gg\doteq\tau(abdc)$ is a
bilinear hermitian positive form on $(H\otimes H)\otimes\mathbb C$
and the following properties are satisfied
\begin{enumerate}
\item the mappings $v\rightarrow v\otimes1$ and
$v\rightarrow1\otimes v$ are isometric embeddings, i.e.
$$
(v,v)=\ll v\otimes1,v\otimes1\gg=\ll1\otimes v,1\otimes v\gg
$$
\item $\ll\cdot,\cdot\gg$ is cyclic in the following sense
$$
\ll a\otimes b,c\otimes d\gg=\ll c\otimes a,d\otimes b\gg
$$
\item $\ll\cdot,\cdot\gg$ is self-adjoint in the following sense
$$
\ll a\otimes b,c\otimes d\gg=\overline{\ll b\otimes a,d\otimes c\gg}
$$
\item $\ll\cdot,\cdot\gg$ verifies the following property
$$
\ll a\otimes b,c\otimes d\gg=\ll b\otimes d,a\otimes c\gg
$$
\item The mapping $J_H:(H\otimes H)\otimes\mathbb C\rightarrow (H\otimes H)\otimes\mathbb C$ defined by
setting $J_H(a\otimes b)=b\otimes a$ is an isometric involution,
i.e.
\begin{enumerate}
\item $J_H(J_H(a\otimes b))=a\otimes b$
\item $\ll J_H(a\otimes b),J_H(a\otimes b)\gg=\ll a\otimes b,a\otimes
b\gg$
\end{enumerate}
\end{enumerate}
\end{prop}
In this article we want to consider Hilbert spaces which verify
these five additional properties. Before giving the definition let
us observe that properties 2. and 3. together imply properties 4.
and 5. Indeed
\begin{lem}
Let $(H,(\cdot,\cdot),\ll\cdot,\cdot\gg)$ be a Hilbert space
equipped with a bilinear positive hermitian form $\ll\cdot,\cdot\gg$
on $(H\otimes H)\otimes\mathbb C$. If $\ll\cdot,\cdot\gg$ verifies
2. and 3. then it verifies also 4. and 5.
\begin{proof}
Suppose 2. and 3. are verified. Applying hermitianity, 2. and
hermitianity again we get 4. Indeed
$$
\ll a\otimes b,c\otimes d\gg=\overline{\ll c\otimes d,a\otimes
b\gg}=\overline{\ll a\otimes c,b\otimes d\gg}=
$$
$$
=\ll b\otimes d,a\otimes c\gg
$$
On the other hand, applying 3. and hermitianity, we get 5. Indeed
$$
\ll J_H(a\otimes b),J_H(a\otimes b)\gg=\ll b\otimes a,b\otimes
a\gg=\overline{\ll a\otimes b,a\otimes b\gg}=
$$
$$
=\ll a\otimes b,a\otimes b\gg
$$
\end{proof}
\end{lem}
Notice that we have not used the completeness with respect to
$(\cdot,\cdot)$. Thus we can give the following
\begin{defin}
A cyclic pre-Hilbert space is a quadruple
$(V,(\cdot,\cdot),1,\ll\cdot,\cdot\gg)$, where $(V,(\cdot,\cdot))$
is a real pre-Hilbert space, $1\in V$ is a pointed vector such that
$(1,1)=||1||^2=1$ and $\ll\cdot,\cdot\gg$ is a bilinear
complex-valued, hermitian positive form on $(V\otimes
V)\otimes\mathbb C$ verifying properties 1.,2. and 3. (and,
consequently, 4. and 5.).
\end{defin}
\setlength{\parskip}{1ex plus 0.5ex minus 0.2ex}

\section{Relation with Connes' embedding conjecture}

We have begun studying cyclic spaces motivated by Connes' embedding
conjecture. Before explaining how they are related to each other,
let us briefly recall Connes' embedding conjecture. Let $R$ be the
hyperfinite $II_1$ factor (with unique trace denoted by $\tau$) and
let $\omega\in\beta(\mathbb N)\setminus\mathbb N$ be a free
ultrafilter on the natural number. One can construct the ultrapower
$R^\omega$ in the following way: first consider
$l^\infty(R)=\{(x_n)_n\subseteq R: sup_n||x_n||<\infty\}$; then
consider its ideal $I_\omega=\{(x_n)_n\in
l^\infty(R^\omega):lim_{n\rightarrow\omega}\tau(x_n^*x_n)^{\frac{1}{2}}=0\}$;
finally consider the quotient $R^\omega= l^\infty(R)/I_\omega$. It
turns out to be a non weakly separable $II_1$ factor with trace
$\tau_{R^\omega}(x+I_\omega)=lim_{n\rightarrow\omega}\tau(x_n)$,
where $(x_n)$ is any representative sequence for $x$. Connes'
embedding conjecture states that any $II_1$-factor with separable
predual embeds into $R^\omega$ (\cite{Co}). This conjecture has
become more and more interesting in recent years, since many authors
have found lots of equivalent conditions showing that this
conjecture is linked to several branches of mathematics (like group
theory and metric geometry), besides being transversal to most of
the sub-specializations of Operator Algebras (see \cite{Br},
\cite{Br2}, \cite{Ca-Pa}, \cite{Co-Dy}, \cite{El-Sz}, \cite{Ha-Wi},
\cite{Ki}, \cite{Ne-Th}, \cite{Oz}, \cite{Pe}, \cite{Ra1},
\cite{Ra2}, \cite{Ra3}, \cite{Vo2}, \cite{Vo3} for some reference).
Here is the problem we want to focus

\begin{prob}\label{cyclic}
Does every separable cyclic space embed into some $II_1$-factor with
separable predual?
\end{prob}

We are interested in this problem because an affirmative answer
would imply the existence of an algorithm to check Connes' embedding
conjecture. Indeed
\begin{enumerate}
\item Take a $II_1$-factor with separable predual $M$. If Prob.\ref{cyclic} has affirmative answer, then we could
theoretically enumerate all the inequalities verified by the moments
of order $3$ and $4$ in $M$. They are positive definite polynomials
of degree less than or equal to 4, that are quite easy to
understand, being exactly the inequalities of a cyclic space.
\item Take these polynomials and calculate their own infimum on
positive matrices of order $n$. Let $\varepsilon_n\geq0$ be such an
infimum. Observe that Connes' embedding conjecture is true if and
only if $\varepsilon_n$ converges to $0$, when $n$ goes to infinity.
Indeed Connes' embedding problem has an affirmative answer if and
only if one can approximate the moments of order $3$ and $4$ (see
\cite{Ra3}).
\end{enumerate}
If the Connes embedding conjecture is true then the algorithm is
infinite for at least one $M$. On the other hand, if it finite for
some $M$, that is, it stops after a finite time, then the Connes
embedding conjecture might be true or false. In this case, the
algorithm could be used as a tool for constructing possible
counter-example.

\section{Extension of cyclic vector spaces}

The idea to answer Prob.\ref{cyclic} is the following: suppose we
have an orthonormal basis $\{x_n\}$ for $L^2(M,\tau)$, then we would
have
$$
x_ix_j=\sum_n\alpha_{ij}^nx_n
$$
and thus the first requirement is that an element of $V\otimes V$
should be actually an element of $V\otimes1$. It means that the
first step is to extend the cyclic structure by adjoining elements.
More precisely we have to extend the cyclic structure on $V$ to a
cyclic structure on a space $W$ of the shape $V\oplus\mathbb RY$,
where $Y$ is an indeterminate, in order to reconstruct step by step
the product. We mean that, chosen arbitrarily $y\in(V\otimes
V)\otimes\mathbb C$, $y$ has to be represented as the indeterminate
$Y$, i.e. $y=1\otimes Y=Y\otimes1$. This is why extending the cyclic
scalar product $\ll\cdot,\cdot\gg$ to one over $((V\oplus\mathbb
RY)\otimes(V\oplus\mathbb RY))\otimes\mathbb C$ (where $Y$
represents an arbitrary element in $V\otimes V$) means exactly that
we are extending the scalar product on $V$ to get the fixed product
$y$ of elements in $V$. Such a purpose forces some necessary
assumptions on $y$:
\begin{enumerate}
\item $y$ must be self-adjoint, in the sense that $J_Vy=y$.
\item $y$ must have norm $1$, i.e. $\ll y,y\gg=1$.
\item $y$ is not an element of $V\otimes1$ or $1\otimes V$
(otherwise we would have trivial product). In this case we say that
$y$ is a \emph{non-trivial element in $V\otimes V$}.
\end{enumerate}
Unfortunately we are not able to extend the structure exactly, but
just approximately.
\begin{defin}
Let $V$ be a finite dimensional cyclic vector space with orthonormal
basis $\{x_1,...x_n\}$. An $\varepsilon$-perturbation of the
original scalar product $\ll\cdot,\cdot\gg$ is another scalar
product $\ll\cdot,\cdot\gg_\varepsilon$ such that
$$
|\ll x_i\otimes x_j,x_k\otimes x_l\gg-\ll x_i\otimes x_j,x_k\otimes
x_l\gg_\varepsilon|<\varepsilon\,\,\,\,\,\,\,\,\,\,\,\,\,\,\,\,\,\,\forall
i,j,k,l\in\{1,...n\}
$$
\end{defin}
\begin{prop}\label{extension}
Let $V$ be a finite dimensional cyclic vector space and $y$ a
self-adjoint and non-trivial element in $(V\otimes V)\otimes\mathbb
C$ with norm 1. For every $\varepsilon>0$, there exists an
$\varepsilon$-perturbation of $\ll\cdot,\cdot\gg$ which extends to a
cyclic structure on $W:=V\oplus\mathbb RY$ with the property
$y=1\otimes Y=Y\otimes1$.
\end{prop}
Proof of this proposition is quite technical, so we will divide it
in several steps. Indeed, let $1,x_2...x_n$ an orthonormal basis of
$V$ (1 is the pointed vector on $V$), we need to define the products
$$
\ll Y\otimes x_i,x_j\otimes
x_l\gg_\varepsilon\,\,\,\,\,\,\,\,\,\,\,\,\,\,\,\,\ll Y\otimes
x_i,Y\otimes x_j\gg_\varepsilon\,\,\,\,\,\,\,\,\,\,\,\,\,\,\,\,\ll
Y\otimes Y,x_i\otimes x_j\gg_\varepsilon
$$
$$
\ll Y\otimes x_i,x_j\otimes
Y\gg_\varepsilon\,\,\,\,\,\,\,\,\,\,\,\,\,\,\,\,\ll Y\otimes
x_i,Y\otimes Y\gg_\varepsilon\,\,\,\,\,\,\,\,\,\,\,\,\,\,\,\,\ll
Y\otimes Y,Y\otimes Y\gg_\varepsilon
$$
The remaining products $\ll x_i\otimes x_j,x_k\otimes
x_l\gg_\varepsilon$ will be defined in the course of the proof, when
we find the suitable $\varepsilon$-perturbation of
$\ll\cdot,\cdot\gg$. The most technical part of the proof is the
definition of $\ll Y\otimes x_i,x_j\otimes x_k\gg_\varepsilon$,
which will be the first and second step. In the first step we follow
a sequence of necessary conditions in order to construct a linear
system whose solutions allow us to define such products; in the
second step we solve this linear system. Before going into the first
step, let us state some
preliminary notions.\\\\
By the fifth property in Prop.\ref{properties}, $J_V$ behaves on
$(V\otimes V)\otimes\mathbb C$ like an involution, so it is natural
to fix the following terminology.
\begin{defin}
An element $x\in (V\otimes V)\otimes\mathbb C$ is called
self-adjoint if $J_Vx=x$. For an element $x\in (V\otimes
V)\otimes\mathbb C$ which is not self-adjoint, its real part is
$Re(x)=\frac{x+J_Vx}{2}$ and the imaginary part is
$Im(x)=\frac{x-J_Vx}{2i}$.
\end{defin}
Step 1\\
Let $P$ be the projection of $(V\otimes V)\otimes\mathbb C$ onto the
first $(V\otimes\mathbb C)^\perp$. Observe that
$$
\ll Y\otimes x_i,x_j\otimes x_k\gg_\varepsilon=\ll Y\otimes
x_i,P(x_j\otimes x_k)+(1-P)(x_j\otimes x_k)\gg_\varepsilon=
$$
$$
=\ll Y\otimes x_i,P(x_j\otimes x_k)\gg_\varepsilon+\ll Y\otimes
x_i,(1-P)(x_j\otimes x_k)\gg_\varepsilon
$$
Consider the second summand
$$
\lambda_{ij}^k=\ll Y\otimes x_i,(1-P)(x_j\otimes
x_k)\gg_\varepsilon=\ll Y\otimes x_i,1\otimes
x_k\gg_\varepsilon=\ll1\otimes Y,x_k\otimes x_i\gg_W
$$
$$
=\ll y,x_k\otimes x_i\gg_\varepsilon=\ll y,x_k\otimes x_i\gg
$$
So we can think of the numbers $\lambda_{ij}^k$ as being
pre-determined. Let us focus on the first summand: we are going to
find a suitable perturbation in order to determine those numbers.
Let $\xi_i$ be the projection of the vector $Y\otimes x_i$ on
$(V\otimes V)\otimes\mathbb C$ and $\eta_i=P\xi_i$. Then
$$
\ll Y\otimes x_i,P(x_j\otimes
x_k)\gg_\varepsilon=\ll\xi_i,P(x_j\otimes
x_k)\gg_\varepsilon=\ll\eta_i,x_j\otimes x_k\gg_\varepsilon
$$
So we would solve our problem if we found suitable $\eta_i$'s. Now
observe that they should verify the following
$$
\ll\eta_i,x_j\otimes x_k\gg_\varepsilon-\ll
J_\varepsilon\eta_j,x_k\otimes
x_i\gg_\varepsilon=\overline{\lambda_{ji}}^k-\lambda_{ij}^k=:\theta_{ij}^k
$$
This is just a linear system. Before attempting to solve it, let us
write it separately for the real and the imaginary part. We get
$$
Re(\theta_{ij}^k)=Re\ll\eta_i,x_j\otimes x_k\gg_\varepsilon-Re\ll
J_\varepsilon\eta_j,x_k\otimes x_i\gg_\varepsilon=
$$
$$
=\frac{1}{2}(\ll\eta_i,x_j\otimes
x_k\gg_\varepsilon+\overline{\ll\eta_i,x_j\otimes
x_k\gg_\varepsilon}-
$$
$$
-\ll J_\varepsilon\eta_j,x_k\otimes x_i\gg_\varepsilon-\overline{\ll
J_\varepsilon x_j,x_k\otimes x_i\gg_\varepsilon})=
$$
$$
=\frac{1}{2}(\ll\eta_i,x_j\otimes x_k\gg_\varepsilon+\ll
J_\varepsilon\eta_i,x_k\otimes x_j\gg_\varepsilon-
$$ $$-\ll J_\varepsilon\eta_j,x_k\otimes
x_i\gg_\varepsilon-\ll\eta_j,x_i\otimes x_k\gg_\varepsilon)=
$$
$$
=\ll\frac{\eta_i+J_\varepsilon\eta_i}{2},\frac{x_j\otimes
x_k+x_k\otimes
x_j}{2}\gg_\varepsilon-\ll\frac{\eta_j+J_\varepsilon\eta_j}{2},\frac{x_k\otimes
x_i+x_i\otimes x_k}{2}\gg_\varepsilon+
$$
$$
-\ll\frac{\eta_i-J_\varepsilon\eta_i}{2i},\frac{x_j\otimes
x_k-x_k\otimes
x_j}{2i}\gg_\varepsilon-\ll\frac{\eta_j-J_\varepsilon\eta_j}{2i},\frac{x_k\otimes
x_i-x_i\otimes x_k}{2i}\gg_\varepsilon+
$$
$$
=\ll Re(\eta_i),Re(x_j\otimes x_k)\gg_\varepsilon-\ll
Im(\eta_i),Im(x_j\otimes x_k)\gg_\varepsilon+
$$
$$
-\ll Re(\eta_j),Re(x_k\otimes x_i)\gg_\varepsilon-\ll
Im(\eta_j),Im(x_k\otimes x_i\gg_\varepsilon
$$
By an analogous calculation we get
$$
Im(\theta_{ij}^k)=\ll Im(\eta_i),Re(x_j\otimes
x_k)\gg_\varepsilon+\ll Re(\eta_i),Im(x_j\otimes
x_k)\gg_\varepsilon+
$$
$$
+\ll Im(\eta_j),Re(x_k\otimes x_i)\gg_\varepsilon-\ll
Re(\eta_j),Im(x_k\otimes x_i\gg_\varepsilon
$$
Thus we have to solve the equations
$$
<(...Re(\eta_i),Im(\eta_i),...,Re(\eta_j),Im(\eta_j)...),v_{ij}^k>=Re(\theta_{ij}^k)
$$
and
$$
<(...Re(\eta_i),Im(\eta_i),...,Re(\eta_j),Im(\eta_j)...),w_{ij}^k>=Im(\theta_{ij}^k)
$$
being
$$
v_{ij}^k=(0,...Re(x_j\otimes x_k),-Im(x_j\otimes
x_k),...0...,-Re(x_k\otimes x_i),-Im(x_k\otimes x_i),...)
$$
and
$$
w_{ij}^k=(0,...Im(x_j\otimes x_k), Re(x_j\otimes x_k),...0...
-Im(x_k\otimes x_i),Re(x_k\otimes x_i)...)
$$
where the non-zero components are exactly the ones
corresponding to $i$ and $j$.\\
Now observe that $i,j$ are switchable everywhere and the case $i=j$
is trivial. So we can suppose $i<j$. Moreover, since the solvability
of a linear system neither depend on permutations of the columns nor
on multiplication by non-zero numbers, we can replace $v_{ij}^k$ and
$w_{ij}^k$ by the following
$$
v_{ij}^k=(0,... Re(x_j\otimes x_k),Im(x_j\otimes
x_k),...0...-Re(x_k\otimes x_i), Im(x_k\otimes x_i),...)
$$
$$
w_{ij}^k=(0,...-Im(x_j\otimes x_k),Re(x_j\otimes
x_k),...0...Im(x_k\otimes x_i),Re(x_k\otimes x_i),...)
$$
Such a re-writing concludes the first step.\\\\
Step 2.\\
The purpose of this step is to find a deformation of the $x_i$'s
(namely: a perturbation of the scalar product) such that the new
$v_{ij}^k$'s, $w_{ij}^k$'s become linearly
independent so that we can solve the equations in Step 1.\\
Suppose we have a linear combination which gives 0:
$$
\sum_{2\leq i<j\leq n,2\leq k\leq
n}\alpha_{ij}^kv_{ij}^k+\sum_{2\leq i<j\leq n,2\leq k\leq
n}\beta_{ij}^kw_{ij}^k=0
$$
Now fix $i$ and look at this relation in the $i$-th component. We
have of course the case $i<j$, but also a contribution that can be
obtained from some $j'<i$. So we can split the previous condition in
the following ones:
$$
\sum_{2\leq i<j\leq n,2\leq k\leq n}(\alpha_{ij}^kRe(x_j\otimes
x_k)-\beta_{ij}^kIm(x_j\otimes x_k))+ $$ $$ +\sum_{2\leq j'<i\leq
n,2\leq k'\leq n}(-\alpha_{j'i}^{k'}Re(x_{k'}\otimes
x_{j'})+\beta_{j'i}^{k'}Im(x_{k'}\otimes x_{j'}))
$$
and
$$
\sum_{2\leq i<j\leq n,2\leq k\leq n}(\alpha_{ij}^kIm(x_j\otimes
x_k)+\beta_{ij}^kRe(x_j\otimes x_k))+ $$ $$ +\sum_{2\leq j'<i\leq
n,2\leq k'\leq n}(\alpha_{j'i}^{k'}Im(x_{k'}\otimes
x_{j'})+\beta_{j'i}^{k'}Re(x_{k'}\otimes x_{j'}))
$$

Now let $s_i$ semicircular (see \cite{Vo}) and $\varepsilon'>0$
small enough. Semicircularity guarantees that
$\sqrt{1-\varepsilon'}x_i\oplus\sqrt{\varepsilon'}s_i=\tilde x_i$
are still an orthonormal basis, for any $\varepsilon'>0$. The choice
of $\varepsilon'$ small enough guarantees that the scalar product
$$
\ll x_i\otimes x_j,x_k\otimes x_l\gg_\varepsilon:=\ll\tilde
x_i\otimes\tilde x_j,\tilde x_k\otimes\tilde x_l\gg
$$
is an $\varepsilon$-deformation. Moreover observe that in this
deformation $x_i\otimes x_k$ are linearly independent and
independent from $x_i\otimes1$. In particular $Re(x_i\otimes x_j)$
and $Im(x_i\otimes x_j)$ are linearly independent over the real
numbers. It follows that we can separate real and imaginary part in
the previous conditions and get
$$
\sum_{2\leq i<j\leq n,k=2,...n}\alpha_{ij}^kRe(x_j\otimes
x_k)-\sum_{2\leq j'<i\leq
n,k'=2,...n}\alpha_{j'i}^{k'}Re(x_{k'}\otimes x_{j'})=0
$$
$$
\sum_{2\leq i<j\leq n,k=2,...n}\alpha_{ij}^kIm(x_j\otimes
x_k)+\sum_{2\leq j'<i\leq
n,k'=2,...n}\alpha_{j'i}^{k'}Im(x_{k'}\otimes x_{j'})=0
$$
$$
-\sum_{2\leq i<j\leq n,k=2,...n}\beta_{ij}^kIm(x_j\otimes
x_k)+\sum_{2\leq j'<i\leq
n,k'=2,...n}\beta_{j'i}^{k'}Im(x_{k'}\otimes x_{j'})=0
$$
$$
\sum_{2\leq i<j\leq n,k=2,...n}\beta_{ij}^kRe(x_j\otimes
x_k)+\sum_{2\leq j'<i\leq
n,k'=2,...n}\beta_{j'i}^{k'}Re(x_{k'}\otimes x_{j'})=0
$$
Now, let us consider the first two conditions. If in the first sum
$i<k$ or in the second sum $i>k'$, the respective terms cannot
cancel each other, so their coefficients must be zero. So one can
have a term in  the first sum equal to one in the second sum only in
case $i>k$, $i<k'$, $k$ corresponds to $j'$ in the second sum and
$k'$ corresponds to $j$ in the first sum. In this case one has
$\alpha_{ij}^k-\alpha_{j'i}^{k'}=0$ from the first condition and
$\alpha_{ij}^k+\alpha_{j'i}^{k'}=0$ from the second one. It follows
that these coefficients must be zero. Similarly we obtain that the
$\beta$'s are equal to zero.\\\\
Step 3.\\
Here we want to define the scalar product
$\ll\cdot,\cdot\gg_\varepsilon$ whenever $Y$ appears twice.
Recalling that the following properties have to be satisfied
\begin{enumerate}
\item $\ll Y\otimes x_i,Y\otimes x_j\gg_\varepsilon=\overline{\ll x_i\otimes Y, x_j\otimes
Y\gg_\varepsilon}$
\item $\ll Y\otimes x_i,Y\otimes x_j\gg_\varepsilon=\overline{\ll Y\otimes Y, x_i\otimes
x_j\gg_\varepsilon}$
\end{enumerate}
it follows that it will be enough to define the numbers $\ll
Y\otimes x_i,Y\otimes x_j\gg_\varepsilon$ and $\ll Y\otimes
x_i,x_j\otimes Y\gg_\varepsilon$. So, we can define the matrix $(\ll
Y\otimes x_i,Y\otimes x_j\gg_\varepsilon)$ as any positive matrix
(Indeed the perturbation in the second step causes $x_i\otimes x_j$
to be linearly independent with respect to
$\ll\cdot,\cdot\gg_\varepsilon$ and then the second of the previous
conditions gives no further constrictions). Finally we can
set $\ll Y\otimes x_i,x_j\otimes Y\gg_\varepsilon=0$.\\\\
Forth Step:\\
We can complete the proof very easily. Indeed we can set $\ll
Y\otimes x_i,Y\otimes Y\gg_\varepsilon=0$, without contradictions.
Finally, Bessel's inequality forces
$$
\ll Y\otimes Y,Y\otimes Y\gg_\varepsilon\geq\sum_{i,j}|\ll Y\otimes
Y,x_i\otimes x_j\gg_\varepsilon|^2
$$
Also in this case there are no contradictions: it is enough to
choose $\ll Y\otimes Y,Y\otimes Y\gg$ large enough. (What is the
smallest possible value?)

\begin{cor}
Let $V$ be a cyclic finite dimensional space with orthonormal basis
$x_1,... x_n$. Then for every $\varepsilon>0$ there exists a
countably generated cyclic space $W_{\varepsilon}$, with cyclic
structure $\ll\cdot,\cdot\gg_\varepsilon$, that verifies the
following properties
\begin{enumerate}
\item $W$ extends $V$ as a vector space, i.e. the set
$\{x_1,...x_n\}$ extends to a basis $\{x_1,...x_n,x_{n+1},...\}$ of
$W$.
\item $\ll\cdot,\cdot\gg_\varepsilon|_V$ is an
$\varepsilon$-deformation of the cyclic structure on $V$.
\item $d_\varepsilon(x_i\otimes x_j,W_\varepsilon\otimes1)=0$ for every
$i,j\in\mathbb N$, where $d_\varepsilon(x,y)=\sqrt{\ll
x-y,x-y\gg_\varepsilon}$
\end{enumerate}

\begin{proof}
It is enough to iterate the previous lemma, taking $\varepsilon/2^n$
at step $n$.
\end{proof}
\end{cor}

\section{Problems we were not able to solve}
Let $V$ be a separable cyclic vector space with orthonormal basis
$\{x_n\}$. We can think about $x_i$ as the operator on $V$ defined
by setting $x_i(x_j)=x_i\otimes x_j$. Prop.\ref{extension}
guarantees that this operator is well defined by linearity, but the
problem is that it could be unbounded. Indeed
$||x_i||^2=\sum_n\alpha_{ii}^n||x_n||^2$ could be infinite. So we
have several open questions

\begin{enumerate}
\item Is the set of operators obtained in such a way a tracial
algebra or at least as unbounded algebra of operators of type $II$
in the sense of Inoue (see \cite{In})?
\item Can we modify the proof of Prop.\ref{extension} in such a way that we get
bounded operators?
\item What is the relation between this construction and that of
Netzer and Thom (see \cite{Ne-Th}), who seem to obtain similar
objects?
\end{enumerate}

\section{Acknowledgement}
The authors are grateful to Robin Hillier for reading the draft of
the paper and for suggesting a correction.

Valerio CAPRARO - University of Rome ''Tor Vergata'' -
capraro@mat.uniroma2.it or valerio.capraro@virgilio.it\\\\
Florin R\u ADULESCU - University of Rome ''Tor Vergata'' -
radulesc@mat.uniroma2.it

\end{document}